\documentclass{article}

\title{ \bf Leonard pairs and the
Askey-\\Wilson relations}
\author{Paul Terwilliger\footnote{University of Wisconsin at
Madison}\\ Raimundas Vidunas\footnote{University of Antwerp;
Supported by the ESF NOG programme, 
and by the EC TMR project.}\\
}

\usepackage{latexsym}
\usepackage{mathrsfs}
\newtheorem{theorem}{Theorem}[section]
\newtheorem{corollary}[theorem]{Corollary}
\newtheorem{lemma}[theorem]{Lemma}

\newtheorem{definition}[theorem]{Definition}
\newtheorem{remark}[theorem]{Remark}
\newtheorem{example}[theorem]{Example}

\newcommand{\proof}{{\bf Proof. }}
\newcommand{\qed}{\hfill {$\Box$}\\} 

\newcommand{\refi}[1]{{\it (#1)}}

\newcommand{\fK}{\mbox{\bf K}}

\newcommand{\CC}{\mbox{\bf C}}
\newcommand{\inint}[2]{0\le #1\le #2}
\newcommand{\equal}{\!\!\!=\!\!\!}
\newcommand{\theorembreak}{}

\begin{document}

\maketitle

\begin{abstract}
\noindent Let $\fK$ denote a field and let $V$ denote a vector
space over $\fK$ with finite positive dimension.
We consider an ordered pair of linear transformations
  $A:V\to V$ and $A^*:V\to V$
 which satisfy the following two properties:
\begin{enumerate}
\item There exists a basis for $V$ with respect to which the matrix
representing $A$ is irreducible tridiagonal and
the matrix representing $A^*$ is
diagonal.
\item There exists a basis for $V$ with respect to which the matrix
representing $A^*$ is irreducible
tridiagonal and the matrix representing $A$ is
diagonal.
\end{enumerate}
We call such a pair a {\em Leonard pair} on $V$.
Referring to the above Leonard pair,
we show there
exists a sequence of scalars
$\beta,\gamma,\gamma^*,\varrho,\varrho^*,
\omega,\eta,\eta^*$ taken from $\fK$ such that
both
\begin{eqnarray*}
A^2 A^*-\beta A A^*\!A+A^*\!A^2-\gamma\left( A A^*\!+\!A^*\!A
\right)-\varrho\,A^* &\equal& \gamma^*\!A^2+\omega A+\eta\,I,\\
A^*{}^2\!A-\beta A^*\!AA^*\!+AA^*{}^2-\gamma^*\!\left(A^*\!A\!+\!A
A^*\right)-\varrho^*\!A &\equal&\gamma A^*{}^2+\omega
A^*\!+\eta^*\!I.
\end{eqnarray*}
The sequence is uniquely determined by the Leonard pair
provided the dimension of $V$ is at least $4$. The equations above are
called the {\em Askey-Wilson relations}.
\end{abstract}

\section{Introduction}

We begin by recalling the notion of a {\it Leonard pair}
\cite{TD00,
shape,
LS99,
  qSerre,
   LS24,
   conform,
    lsint,
Terint,
TLT:split,
 TLT:array,
qrac}.
We will use the following terms. Let $X$ denote a square matrix.
Then $X$ is called {\it tridiagonal} whenever each nonzero
entry lies on either the diagonal, the subdiagonal, or the
superdiagonal. Assume $X$ is tridiagonal. Then $X$ is
called {\it irreducible} whenever each entry on
the subdiagonal is nonzero and each entry on the superdiagonal
is nonzero.

\medskip
\noindent We now define a Leonard pair.
For the rest of this paper $\fK$ denotes a field.
\begin{definition} \rm
\cite{LS99}
\label{deflp}
Let $V$ denote a vector space over $\fK$ with finite positive
dimension. By a {\em Leonard pair} on $V$ we mean an ordered pair
$(A,A^*)$, where $A:V\to V$ and $A^*:V\to V$ are linear
transformations which satisfy the following two properties:
\begin{enumerate}
\item There exists a basis for $V$ with respect to which the matrix
representing $A$ is irreducible tridiagonal
and the matrix representing $A^*$ is
diagonal.
\item There exists a basis for $V$ with respect to which the matrix
representing $A^*$ is irreducible tridiagonal
and the matrix representing $A$ is
diagonal.
\end{enumerate}
\end{definition}

\begin{remark} \rm
According to a common notational convention, $A^*$ denotes
 the conjugate transpose of $A$.
We are not using this convention.
 In a Leonard pair $(A,A^*)$ the
linear transformations $A$ and $A^*$ are arbitrary subject to the
conditions {\it (i)} and {\it (ii)} above.
\end{remark}

\noindent
Leonard pairs occur in
combinatorics
\cite{
Cau,
CurNom,
TD00,
TersubI,
lsint},
 representation theory
\cite{
GYZnature,
Zhed,
  GYZlinear,
Koelink3,
cite37,
cite38,
cite39,
cite40,
Rosengren,
lsint},
 and the theory of orthogonal
polynomials \cite{
BanIto,
Leon,
lsint,
TLT:array,
qrac}.
The connection to polynomials is as
follows: there is a natural correspondence between Leonard pairs
and a class of orthogonal polynomials consisting of the $q$-Racah polynomials
\cite{AWil, GR}
and some related polynomials of the Askey-scheme
\cite{KoeSwa,TLT:array,qrac}.
This correspondence is illustrated in the following example.

\begin{example} \label{krawtchouk} \rm
Let $d$ denote a nonegative integer and
define a set $\Omega=\{0,1,2,\dots,d\}$.
Let $V$ denote the vector space  over the complex field $\CC$
consisting
of all functions from $\Omega$ to $\CC$.
The cardinality of $\Omega $ is $d+1$ so $\dim V=d+1$.
Consider the Krawtchouk polynomials \cite[Chapter
1.10]{KoeSwa}:
\begin{equation}
y_n(x)=K_n(x;p,d)={}_2F_1\left( {-n,\,-x\, \atop -d}
;\,\frac{1}{p}\, \right), \qquad\quad 0\le n\le d.
\end{equation}
We view $y_0,y_1,\ldots,y_d$
as elements of $V$.
The Krawchouk polynomials satisfy the following three-term recurrence.
For $0 \leq n \leq d$ and for $x \in \Omega$,
\begin{equation} \label{krawrec}
-x\,y_n(x)=p\,(d\!-\!n)\,y_{n+1}(x)-\big(p\,(d\!-\!n)+(1\!-\!p)\,n
\big)\,y_n(x)+(1\!-\!p)\,n\,y_{n-1}(x),
\end{equation}
where $y_{-1}(x)=0$ and $y_{d+1}(x)=0$. The
Krawtchouk polynomials satisfy the following difference equation.
For $0 \leq n \leq d$ and for $x \in \Omega$,
\begin{equation} \label{krawdif}
-n\,y_n(x)=p\,(d\!-\!x)\,y_n(x\!+\!1)-\big(p\,(d\!-\!x)+(1\!-\!p)\,x
\big)\,y_n(x)+(1\!-\!p)\,x\,y_n(x\!-\!1),
\end{equation}
where $y_n(-1)=0$ and $y_n(d+1)=0$.
We now define
linear transformations $A:V\to V$ and $A^*:V\to V$.
We begin with $A$.
 For
$f\in V$, $Af$ is the element of $V$ which satisfies
\[
(Af)(x) = x\,f(x),  \qquad  \qquad x \in \Omega.
\]
We now define $A^*$. For $f \in V$,
$A^*f$ is the element in $V$ which satisfies
\[
(A^*f)(x)=p\,(d\!-\!x)\,f(x\!+\!1)-\big(
p\,(d\!-\!x)+(1\!-\!p)\,x\big)\,f(x)+
(1\!-\!p)\,x\,f(x\!-\!1)
\]
for $x \in \Omega$.
We mentioned the polynomials $y_0,y_1,\ldots,y_d$
are elements of $V$. These elements form a basis of $V$ since they
are linearly independent. According to (\ref{krawrec}), the matrix
representing $A$ in this basis is irreducible tridiagonal:
\[
\left( \begin{array}{ccccc} pd & p-1 & & &{\mathbf 0}\\
-pd & p(d\!-\!2)\!+\!1 & 2(p-1) \\
& -p(d\!-\!1) & p(d\!-\!4)\!+\!2 & \ddots \\
& & \ddots & \ddots & d(p-1) \\
{\mathbf 0}& & & -p & -pd+d
\end{array} \right).
\]
According to (\ref{krawdif}), the matrix representing $A^*$ in
this basis is
$\mbox{diag}(0,-1,-2,\ldots,-d)$.
For $0\le n\le d$ let $y^*_n$ denote the element of $V$ which
satisfies
\[
y^*_n(x)=\delta_{nx}, \qquad \qquad \mbox{$x \in \Omega$,}
\]
where $\delta_{nx}$ is the Kronecker delta. The sequence
$y^*_0,y^*_1,\ldots,y^*_d$ forms a basis of $V$.
With respect to this basis the matrix representing
$A$ is
$\mbox{diag}(0,1,2,\ldots,d)$ and the matrix representing $A^*$
is irreducible tridiagonal:
\[
\left( \begin{array}{ccccc} -pd & 1-p  & & & {\mathbf 0}\\
pd & p(2\!-\!d)\!-\!1 & 2(1-p) \\
& p(d\!-\!1) & p(4\!-\!d)\!-\!2 & \ddots \\
& & \ddots & \ddots & d(1-p) \\
{\mathbf 0}& & & p & pd-d
\end{array} \right).
\]
Therefore $(A,A^*)$ is
a Leonard pair on $V$. \qed
\end{example}

\noindent
In order to motivate our main result we cite a theorem of
Terwilliger.

\begin{theorem}
 {\rm \cite[Theorem 1.12]{LS99}}
\label{tdptheorem}
Let $V$ denote a vector space over $\fK$ with finite positive
dimension. Let $(A,A^*)$ denote
a Leonard pair on $V$. Then there
exists a sequence of scalars
$\beta,\gamma,\gamma^*,\varrho,\varrho^*$ taken from $\fK$ such
that both
\begin{eqnarray} \label{askwilb1}
\lbrack\, A,\, A^2 A^*-\beta A A^*\!A+A^*\! A^2-\gamma\left(
A A^*\!+\!A^*\!A\right)-\varrho\,A^* \,\rbrack & = & 0,\\
\label{askwilb2} \lbrack A^*\!,\, A^*{}^2\!A-\beta A^*\!AA^*
\!+AA^*{}^2-\gamma^*\!\left(A^*\!A\!+\!AA^*\right)
-\varrho^*\!A\,\rbrack&=&0.
\end{eqnarray}
Here $[r,s]$ means $rs-sr$. The sequence is uniquely
determined by the pair $(A,A^*)$ provided the dimension of $V$ is
at least $4$.
\end{theorem}

\noindent
The equations (\ref{askwilb1}), (\ref{askwilb2}) are called the {\em
tridiagonal relations} \cite[Lemma 5.4]{TersubIII}.
See
\cite{
TD00,
TersubI,
TersubII,
LS99,
qSerre}
for more information on these relations.

\medskip
\noindent The main result of this paper is the following extension
of Theorem
\ref{tdptheorem}.

\begin{theorem} \label{lptheorem}
Let $V$ denote a vector space over $\fK$ with finite positive
dimension. Let $(A,A^*)$ denote
a Leonard pair on $V$. Then there
exists a sequence of scalars
$\beta,\gamma,\gamma^*,\varrho,\varrho^*$, $\omega,\eta,\eta^*$
taken from $\fK$ such that both
\begin{eqnarray}  \label{askwil1}
A^2 A^*-\beta A A^*\!A+A^*\!A^2-\gamma\left( A A^*\!+\!A^*\!A
\right)-\varrho\,A^* &\equal& \gamma^*\!A^2+\omega A+\eta\,I,\\
\label{askwil2} A^*{}^2\!A-\beta A^*\!AA^*\!+AA^*{}^2-
\gamma^*\!\left(A^*\!A\!+\!A A^*\right)-\varrho^*\!A &\equal&
\gamma A^*{}^2+\omega A^*\!+\eta^*I.
\end{eqnarray}
The sequence is uniquely determined by the pair $(A,A^*)$ provided
the dimension of $V$ is at least $4$.
\end{theorem}

\noindent
As far as we know,
the relations (\ref{askwil1}), (\ref{askwil2}) first appeared
in
\cite{Zhidd}.
In that article it is shown that the Askey-Wilson
polynomials give a pair of infinite matrices which
satisfy
 (\ref{askwil1}), (\ref{askwil2}).
See
\cite{GYZnature,
Zhed,
  GYZlinear,
ZheCart}
  for related work.
In these articles the relations
 (\ref{askwil1}), (\ref{askwil2}) are called
the {\em
Askey-Wilson relations}  and we shall also use this term.
One of the
relations
 (\ref{askwil1}), (\ref{askwil2})
shows up in work of
Gr{\"u}nbaum
and Haine on the ``bispectral problem''
\cite{GH1}
where it is called a {\it $q$-analog of the  string equation}.
See \cite{GH4, GH5,GH7, GH6, GH3, GH2}
for related work.

\medskip
\noindent The plan for the rest of this paper is as follows.
 We will first give a proof of
  Theorem \ref{tdptheorem} which is considerably shorter
  than the one in \cite{LS99}. We will then display some formulae
  which can be used to compute the scalars from
  Theorem \ref{tdptheorem}.
After this  we will use
  Theorem \ref{tdptheorem}
to obtain Theorem \ref{lptheorem}.
We will then display some formulae which can be used
to compute the scalars
from Theorem
 \ref{lptheorem}.
Finally we will illustrate Theorem
\ref{lptheorem} using Example
 \ref{krawtchouk}.

\section{Preliminaries}

\noindent In this section we review some notation and basic concepts.
Let $d$ denote a nonnegative integer.
Let $V$ denote a vector space over $\fK$ with dimension $d+1$.
We let $\mbox{End}(V)$ denote the $\fK$-algebra consisting
of all linear transformations from $V$ to $V$.
For $A \in \mbox{End}(V)$,
by an {\em eigenvalue} of $A$ we mean a root of the characteristic
polynomial of $A$. The eigenvalues of $A$ are contained in
the algebraic closure of $\fK$.
 We say that $A$ is {\em multiplicity-free}
whenever it has $d+1$ distinct eigenvalues, all of which lie
in $\fK$.
Assume $A$ is multiplicity-free. Let
$\theta_0,\theta_1,\ldots,\theta_d$ denote an ordering of the
eigenvalues of $A$, and for $\inint{i}{d}$ put
\[
E_i=\prod_{0\le j\le d \atop j\neq i}
\frac{A-\theta_j\,I}{\theta_i-\theta_j},
\]
where $I$ denotes the identity in
$\mbox{End}(V)$.
By elementary algebra,
\begin{eqnarray} \label{acomei}
A\,E_i=E_i\,A=\theta_i\,E_i, & & \qquad  0\le i\le d,\\
E_i\,E_j=\delta_{ij}\,E_i, & &  \qquad 0\le i,j\le d,
\\
\label{sumei} \sum_{i=0}^d E_i=I, & & A=\sum_{i=0}^d \theta_i\,E_i.
\end{eqnarray}
Let $\mathcal D$ denote the $\fK$-subalgebra of $\mbox{End}(V)$
generated by  $A$. Using
(\ref{acomei})--(\ref{sumei})
we find
$E_0,E_1,\ldots,E_d$ form a
basis for the $\fK$-vector space $\mathcal D$.
 We refer to
$E_i$ as the {\em primitive idempotent} of $A$ associated with
$\theta_i$. We have
\begin{equation}
V=E_0 V+E_1 V+\ldots+E_d V \qquad \qquad (\mbox{direct sum}).
\end{equation}
For $\inint{i}{d}$, $E_i V$ is the (one-dimensional)
eigenspace of
$A$  associated with the eigenvalue $\theta_i$, and $E_i$
acts on $V$ as the projection onto this eigenspace.
We remark that $\lbrace A^i|0 \leq i \leq d\rbrace $
is a basis for the $\fK$-vector space $\mathcal D$,
and that
$\prod_{i=0}^d (A-\theta_iI)=0$.

\medskip
\noindent Later in this paper we will encounter
sequences of scalars which satisfy a certain  recurrence.
We take a moment to discuss this recurrence.
Let $d$ denote a nonnegative integer and let
$\theta_0, \theta_1, \ldots, \theta_d$ denote a sequence of
 scalars taken from $\fK$.
Given $\beta \in \fK$, we say
the sequence
$\theta_0, \theta_1, \ldots, \theta_d$
is {\it $\beta$-recurrent}
whenever
\begin{eqnarray}
\label{eq:1rec}
\theta_{i-2}
-(\beta +1)
\theta_{i-1}
+(\beta +1)
\theta_i
-\theta_{i+1}=0
\end{eqnarray}
for $2 \leq i \leq d-1$.
Given scalars $\beta, \gamma$ in $\fK$,
we say
the sequence
$\theta_0, \theta_1, \ldots, \theta_d$
is {\it $(\beta, \gamma)$-recurrent}
whenever
\begin{eqnarray}
\label{eq:2rec}
 \theta_{i-1}-\beta\theta_i +\theta_{i+1} = \gamma
\end{eqnarray}
for $1 \leq i \leq d-1$. Observe that for  $\beta \in \fK$ the
following are equivalent:
 {\it (i)} the sequence $\theta_0,
\theta_1, \ldots, \theta_d$ is $\beta$-recurrent;
 {\it (ii)} there
exists $\gamma \in \fK$ such that $\theta_0, \theta_1, \ldots,
\theta_d$ is $(\beta, \gamma)$-recurrent. We also have the
following.
\begin{lemma}
\label{lem:pre}
Let $d$ denote a nonnegative integer and let
$\theta_0, \theta_1, \ldots, \theta_d$ denote a sequence of
 scalars taken from $\fK$.
Let $\beta, \gamma$ denote scalars in $\fK$.
Then the following (i), (ii) hold.
\begin{enumerate}
\item Assume
 $\theta_0, \theta_1, \ldots, \theta_d$
is $(\beta, \gamma)$-recurrent.
Then there  exists $\varrho \in \fK$ such that
\begin{eqnarray}
\label{eq:3form}
\theta^2_{i-1}-\beta \theta_{i-1}\theta_i+\theta^2_i
-\gamma (\theta_{i-1}+\theta_i)=\varrho,
\qquad 1 \leq i \leq d.
\end{eqnarray}
\item Assume there exists
 $\varrho \in \fK$ such that
{\rm (\ref{eq:3form})} holds. Further assume
$\theta_{i-1}\not=\theta_{i+1} $ for $1 \leq i \leq d-1$. Then the
sequence
  $\theta_0, \theta_1, \ldots, \theta_d$
is $(\beta, \gamma)$-recurrent.
\end{enumerate}
\end{lemma}
\proof Let $p_i$ denote the
left-hand side of (\ref{eq:3form}) and observe
\begin{eqnarray*}
p_{i} - p_{i+1} =(\theta_{i-1}-\theta_{i+1})
(\theta_{i-1}-\beta \theta_i + \theta_{i+1}-\gamma)
\end{eqnarray*}
for $1 \leq i \leq d-1$.
Our assertions {\it (i), (ii)}
are  routine consequences of this.
\qed

\section{General setting}

In this section we establish some basic results concerning
Leonard pairs.
We begin with a comment.

\begin{lemma}
{\rm \cite[Lemma 1.3]{LS99}}
 \label{multfree} Let $V$ denote a vector
space over $\fK$ with finite positive dimension. Let $(A,A^*)$
denote a Leonard pair on $V$. Then each of $A$, $A^*$ is
multiplicity-free.
\end{lemma}
\proof Set $d+1=\dim V$. Recall that there exists a basis
for $V$ with respect to which the matrix representing $A$ is
diagonal. Hence the eigenvalues of $A$ are in $\fK$. Moreover the
degree of the minimal polynomial of $A$
is equal to the number of distinct
eigenvalues of $A$. Recall there exists a basis for $V$
with respect to which the matrix representing $A$ is irreducible
tridiagonal. From the shape of this matrix we find that
$I,A,
A^2,\ldots,A^d$ are linearly independent. Therefore the degree of
the minimal polynomial of $A$ is equal to
$d+1$. It follows that the
eigenvalues for $A$ are mutually distinct. We have  now
shown that $A$
is multiplicity-free. 
Applying this argument to the Leonard pair $(A^*,A)$ we
find that $A^*$ is multiplicity-free. \qed

\noindent
For the rest of this paper we adopt the following
notational convention.

\begin{definition} \label{primdef} \rm
Let $d$ denote a nonnegative integer and
let $V$ denote a vector
space over $\fK$ with dimension $d+1$. Let $(A,A^*)$ denote a Leonard
pair on $V$. Let $v_0,v_1,\ldots,v_d$ denote a basis for $V$ which
satisfies condition \refi{ii} of Definition \ref{deflp}. For
$\inint{i}{d}$ the vector $v_i$ is an eigenvector of $A$; let
$\theta_i$ (resp. $E_i$) denote the corresponding eigenvalue
(resp. primitive idempotent).
Let $v^*_0,v^*_1,\ldots,v^*_d$ denote a basis for $V$
which satisfies condition \refi{i} of Definition \ref{deflp}. For
$\inint{i}{d}$ the vector $v^*_i$ is an eigenvector of $A^*$; let
$\theta^*_i$ (resp. $E^*_i$) denote the corresponding eigenvalue
(resp. primitive idempotent).
Let the sequence
$a_0, a_1, \ldots, a_d$ denote  the diagonal of
the matrix which represents $A$ with respect to
 $v^*_0,v^*_1,\ldots,v^*_d$.
Let the sequence
$a^*_0, a^*_1, \ldots, a^*_d$ denote the diagonal of
the matrix which represents $A^*$ with respect to
 $v_0,v_1,\ldots,v_d$.
We remark $a_i = \mbox{tr}(E^*_iA)$ and $a^*_i=\mbox{tr}(E_iA^*)$ for
$0 \leq i \leq d$.
\end{definition}

\begin{lemma} \label{leosystem}
With reference to Definition
 $\ref{primdef}$, the following (i), (ii) hold.
\begin{enumerate}
\item $E^*_i\,A\,E^*_j=\left\{ \begin{array}{cl}
 0, & \mbox{if } |i-j|>1;\\ \neq 0, &  \mbox{if } |i-j|=1,
 \end{array} \right.\qquad 0\le i\le d$.
\item $E_i\,A^*\,E_j=\left\{ \begin{array}{cl}
 0, & \mbox{if } |i-j|>1;\\ \neq 0, &  \mbox{if } |i-j|=1,
 \end{array} \right.\qquad 0\le i\le d$.
\end{enumerate}
\end{lemma}
\proof Assertion \refi{i} follows from the irreducible tridiagonal
shape of the matrix representing $A$ in the basis
$v^*_0,v^*_1,\ldots,v^*_d$. Assertion \refi{ii} is similarly obtained.\qed

\begin{lemma}
\label{lem:ais}
With reference to
Definition
 $\ref{primdef}$, the following (i), (ii) hold.
\begin{enumerate}
\item $E^*_iAE^*_i=a_iE^*_i, \qquad 0 \leq i\leq d.$
\item $E_iA^*E_i=a^*_iE_i, \qquad 0 \leq i\leq d.$
\end{enumerate}
\end{lemma}
\proof
{\it (i)} Abbreviate ${\mathscr A}=\mbox{End}(V)$.
Observe
$E^*_i {\mathscr A}E^*_i$
is a 1-dimensional subspace of the $\fK$-vector space
 $\mathscr A$. The element $E^*_i$ is  nonzero
and contained in
$E^*_i {\mathscr A}E^*_i$
so it spans
$E^*_i {\mathscr A}E^*_i$.
Observe $E^*_iAE^*_i \in
E^*_i {\mathscr A}E^*_i$, so there exists
$\alpha_i \in \fK$ such that
 $E^*_iAE^*_i=\alpha_i E^*_i$.
In this equation we take the trace of both sides
and use $\mbox{tr}(XY)=
 \mbox{tr}(YX)$ to obtain $\alpha_i=a_i$.
The result follows.
\\
\noindent {\it (ii)} Similar to the proof of {\it (i)} above.
\qed

\noindent The following lemma gives some consequences of Lemma
\ref{leosystem}{\it (i)} which we will find useful. Of course
Lemma \ref{leosystem}{\it (ii)}
 has similar consequences.
\begin{lemma} \label{eiarej}
With reference to Definition
 $\ref{primdef}$, the following (i)--(iii) hold for $0 \leq i,j\leq d$.
\begin{enumerate}
\item $E^*_i\,A^r\,E^*_j=0 \qquad \mbox{for} \quad 0 \leq r<|i-j|$.
\item $E^*_i\,A^r\,E^*_j \not=0  \qquad \mbox{for} \quad r=|i-j|$.
\item For $0 \leq r,s\leq d$,
\begin{eqnarray*}
E^*_i A^r A^* A^s E^*_j=\left\{ \begin{array}{cl}
 \theta^*_{j+s}\,E^*_i A^{r+s} E^*_j,&  \mbox{if } i-j=r+s,\\
 \theta^*_{j-s}\,E^*_i A^{r+s} E^*_j,&  \mbox{if } j-i=r+s,\\
 0,& \mbox{if } |i-j|>r+s.\end{array} \right.
\end{eqnarray*}
\end{enumerate}
\end{lemma}
\proof These are routine consequences of the irreducible
tridiagonal shape of the matrix $A$ in the basis
$v^*_0, v^*_1, \ldots, v^*_d$. \qed
%

\section{The proof of Theorem 1.4}
\label{mainproof}

In this section we establish
 Theorem \ref{tdptheorem}.
 We start in a
fashion similar to \cite{LS99}.
\begin{lemma}
{\rm \cite[Lemma 12.1]{LS99}}
 \label{paul12p1} With reference to
Definition $\ref{primdef}$,
let $\mathcal D$ denote the $\fK$-subalgebra of
$\mbox{End}(V)$ generated by $A$.
Then
\begin{equation} \label{span12p1}
\mbox{\rm Span}\{ XA^*Y-YA^*X \,|\, X,Y\in{\cal D} \}= \{
ZA^*-A^*Z \,|\, Z\in{\cal D} \}.
\end{equation}
\end{lemma}
\proof For notational convenience set $E_{-1}=0$ and $E_{d+1}=0$.
We claim that for $0\le i\le d$,
\begin{equation} \label{eiei12p1}
E_iA^*E_{i+1}-E_{i+1}A^*E_i = L_iA^*-A^*L_i,
\end{equation}
where $L_i=E_0+E_1+\ldots+E_i$. To see this, observe by
(\ref{sumei}) and Lemma \ref{leosystem}{\it (ii)} that for $0\le
j\le d$ both
\begin{eqnarray}
E_jA^* & = & E_jA^*E_{j-1}+E_jA^*E_j+E_jA^*E_{j+1} \label{ejas12p1} \\
A^*E_j & = & E_{j-1}A^*E_j+E_jA^*E_j+E_{j+1}A^*E_j.
\label{asej12p1}
\end{eqnarray}
Summing both (\ref{ejas12p1}) and (\ref{asej12p1}) over
$j=0,1,\ldots,i$, and taking the difference between these two sums
gives (\ref{eiei12p1}).
We mentioned earlier that $E_0, E_1, \ldots, E_d$
form a basis for the $\fK$-vector space $\mathcal D$.
From this we find
$L_0,L_1,\ldots,L_d$ form a basis for $\mathcal D$.
We may now argue
\begin{eqnarray*}
&&\hspace{-2cm} \mbox{Span}\{ XA^*Y-YA^*X \,|\, X,Y\in{\cal D} \} \\
&&=\mbox{\rm Span}\{E_iA^*E_j-E_jA^*E_i \,|\, 0\le i,j\le d \} \\
&&=\mbox{\rm Span}\{E_iA^*E_{i+1}-E_{i+1}A^*E_i\,|\,0\le i\le d \} \\
&&=\mbox{\rm Span}\{L_iA^*-A^*L_i \,|\,0\le i\le d \} \\
&&=\{ ZA^*-A^*Z \,|\, Z\in{\cal D} \}
\end{eqnarray*}
and we are done. \qed

\noindent In order to state the next lemma we introduce some
notation.

\begin{definition}
\label{def:poly}
\rm
Given scalars $\beta, \gamma, \varrho$ in $\fK$ we define
a polynomial
\begin{eqnarray*}
P(x,y)&=&x^2-\beta xy+y^2-\gamma(x+y)-\varrho.
\end{eqnarray*}
Given scalars $\beta, \gamma^*, \varrho^*$ in $\fK$ we define
a polynomial
\begin{eqnarray*}
P^*(x,y)&=&x^2-\beta xy+y^2-\gamma^*(x+y)-\varrho^*.
\end{eqnarray*}
\end{definition}

\begin{lemma} \label{paul12p2}
Let $\beta, \gamma, \varrho$ denote scalars in $\fK$.
Then  with reference to
Definition
 $\ref{primdef}$ and
Definition $\ref{def:poly}$,
\begin{eqnarray} \label{askwilb1pre}
\lbrack\, A,\, A^2 A^*-\beta A A^*\!A+A^*\! A^2-\gamma\left(
A A^*\!+\!A^*\!A\right)-\varrho\,A^* \,\rbrack & = & 0
\end{eqnarray}
 if and only if
$P(\theta_{i-1}, \theta_i)=0$ for $1 \leq i \leq d$.
\end{lemma}
\noindent {\bf Proof.}
Let $C$ denote the expression on
the left in
(\ref{askwilb1pre})
and observe
\begin{eqnarray}
\label{eq:csum}
C = \sum_{i=0}^d \sum_{j=0}^d E_iCE_j.
\end{eqnarray}
For $0 \leq i,j\leq d$ we evaluate $E_iCE_j$ using
(\ref{acomei}),
(\ref{askwilb1pre}) and get
\begin{eqnarray}
\label{eq:cpiece}
E_iCE_j = (\theta_i-\theta_j)P(\theta_i, \theta_j)E_iA^*E_j.
\end{eqnarray}
First assume
(\ref{askwilb1pre}) holds,
so that $C=0$.
We show $P(\theta_{i-1},\theta_i)=0$ for
$1 \leq i \leq d$.
Let $i$ be given.
Observe $E_{i-1}CE_i=0$
so
$(\theta_{i-1}-\theta_i)P(\theta_{i-1},\theta_i)E_{i-1}A^*E_i=0$
in view of
(\ref{eq:cpiece}).
Observe
$\theta_{i-1}\not=\theta_i$
by
Lemma
\ref{multfree}
 and
$E_{i-1}A^*E_i\not=0$
by Lemma
 \ref{leosystem}{\it (ii)}, so
$P(\theta_{i-1},\theta_i)=0$.
We have now proved the lemma in one direction.
To obtain the converse,
 assume
$P(\theta_{i-1},\theta_i)=0$ for $1 \leq i \leq d$.
Since $P(x,y)$ is symmetric in $x,y$ we have 
$P(\theta_{i},\theta_{i-1})=0$
for $1 \leq i \leq d$.
We show $C=0$.
By
(\ref{eq:csum}) it suffices to show
$E_iCE_j=0$ for
 $0 \leq i,j\leq d$. Let $i,j$ be given.
We show at least one of the factors on the right in
(\ref{eq:cpiece}) is zero. If $|i-j|>1$ then $E_iA^*E_j=0$ by
Lemma \ref{leosystem}{\em (ii)}. If $|i-j|=1$ then $P(\theta_i,
\theta_j)=0$. If $i=j$ then $\theta_i-\theta_j=0$. We have now
shown at least one of the factors on the right in
(\ref{eq:cpiece}) is zero. 
Therefore $E_iCE_j=0$.
Now each term on the right in
(\ref{eq:csum}) is zero so $C=0$. We now have  (\ref{askwilb1pre}). 
\qed

\begin{corollary} \label{lpscalars1}
Let $\beta,\gamma,\gamma^*,\varrho,\varrho^*$ denote scalars in
$\fK$. Then with reference to
 Definition $\ref{primdef}$
 and Definition $\ref{def:poly}$,
 the following (i), (ii) are equivalent.
\begin{enumerate}
\item The sequence $\beta,\gamma,\gamma^*,\varrho,\varrho^*$
satisfies {\rm (\ref{askwilb1})} and {\rm (\ref{askwilb2})}.
\item for $1 \leq i \leq d$ both
\begin{eqnarray*}
P(\theta_{i-1}, \theta_i)=0, \qquad \qquad
P^*(\theta^*_{i-1}, \theta^*_i)=0.
\end{eqnarray*}
\end{enumerate}
\end{corollary}
\proof
 Apply Lemma
 \ref{paul12p2} to both $(A,A^*)$ and $(A^*,A)$.
\qed

\medskip

\noindent{\bf Proof of Theorem \ref{tdptheorem}.}
 Set $d+1=\dim V$.
 For $\inint{i}{d}$, let $\theta_i$ and $E_i$ (resp.
$\theta^*_i$ and $E^*_i$) denote the eigenvalues and the primitive
idempotents of $A$ (resp. $A^*$) as in Definition \ref{primdef}.
 Let $\mathcal D$ denote the $\fK$-subalgebra of
$\mbox{End}(V)$ generated by $A$.

\theorembreak
 We first assume that $d\ge 3$. By Lemma
\ref{paul12p1} (with $X=A^2$ and $Y=A$) there exists $Z\in\cal D$
such that
\begin{equation} \label{asliebpa}
A^2\,A^*\,A-A\,A^*\,A^2 = Z\,A^*-A^*\,Z.
\end{equation}
Recall $\mathcal D$ has a basis $I,A,A^2,\ldots,A^d$ so
there exists a
 polynomial $p\in\fK[x]$ which has degree at most $d$
and satisfies
$Z=p(A)$. Let $k$ denote the degree of $p$.

\theorembreak
We show $k=3$. We first suppose $k>3$ and obtain a
contradiction. We multiply each term in (\ref{asliebpa}) on the
left by $E^*_k$ and on the right by $E^*_0$. We evaluate the
result using
 (\ref{acomei}) and Lemma
\ref{eiarej}  to find
$0=c\left(\theta^*_0-\theta^*_k\right)E^*_kA^kE^*_0$,
where $c$ denotes the leading coefficient of
$p$. The scalars
$c$ and
$ \theta^*_0-\theta^*_k$ are nonzero
by the construction. Moreover
$E^*_kA^kE^*_0$ is nonzero by
Lemma \ref{eiarej}.
Therefore $ 0\not=c\left(\theta^*_0-\theta^*_k\right)E^*_kA^kE^*_0$
for
a contradiction.
We have now shown $k\leq 3$.
We next assume $k<3$ and obtain a contradiction.
We multiply each term in
(\ref{asliebpa})
on the left by $E^*_3$ and on the right by
$E^*_0$. We evaluate the result using
 (\ref{acomei}) and Lemma
\ref{eiarej} to  find
$\left(\theta^*_1-\theta^*_2\right)E^*_3A^3E^*_0=0$.
The scalar
$ \theta^*_1-\theta^*_2$ is nonzero.
Moreover
$E^*_3A^3E^*_0$ is nonzero by
Lemma \ref{eiarej}.
Therefore
$\left(\theta^*_1-\theta^*_2\right)E^*_3A^3E^*_0\not=0$
for a contradiction.
We have now shown $k=3$.

\theorembreak
We divide both sides of (\ref{asliebpa}) by $c$. The
result is
\begin{equation} \label{prebrack}
(\beta\!+\!1)\!\left(A^2\!A^*\! A\!-\!A A^*\!A^2\right) =
A^3\!A^*\! -\!A^*\!A^3 -\gamma
\left(A^2\!A^*\!-\!A^*\!A^2\right)-\varrho
\left(AA^*\!-\!A^*\!A\right),
\end{equation}
where $\beta=c^{-1}-1$ and where  $\gamma,\varrho$ are appropriate
scalars in $\fK$. From
(\ref{prebrack})
 we routinely obtain (\ref{askwilb1}).
Concerning (\ref{askwilb2}),
pick any integer $i$ $(2 \leq i\leq d-1)$.
We multiply each term in
 (\ref{askwilb1}) on the left by
$E^*_{i-2}$ and on the right by $E^*_{i+1}$. We evaluate
the result using
 (\ref{acomei}) and Lemma
\ref{eiarej} to  find $E^*_{i-2}A^3E^*_{i+1}$ times
\begin{eqnarray}
\label{eq:brecur}
\theta^*_{i-2}-(\beta +1)\theta^*_{i-1}+
(\beta +1)\theta^*_i -\theta^*_{i+1}
\end{eqnarray}
is zero.
Observe $E^*_{i-2}A^3E^*_{i+1}\not=0$ by
Lemma
\ref{eiarej}
so
(\ref{eq:brecur}) is zero.
Apparently the sequence $\theta^*_0, \theta^*_1, \ldots, \theta^*_d$
is $\beta$-recurrent. Therefore there exists
$\gamma^* \in \fK$ such that
$\theta^*_0, \theta^*_1, \ldots, \theta^*_d$
is $(\beta, \gamma^*)$-recurrent.
By Lemma
\ref{lem:pre}
there exists
$\varrho^* \in \fK$ such that
$P^*(\theta^*_{i-1}, \theta^*_i)= 0$
for $1 \leq i \leq d$, where
$P^*$ is from
Definition
\ref{def:poly}.
By this and Lemma \ref{paul12p2} we find $\beta, \gamma^*,
\varrho^*$ satisfy (\ref{askwilb2}).

\theorembreak
We have now shown there exists a sequence of scalars
$\beta, \gamma, \gamma^*, \varrho, \varrho^*$ taken from $\fK$
which satisfies
  (\ref{askwilb1}),
  (\ref{askwilb2}).
We show this sequence is unique.
Let
 $\beta, \gamma, \gamma^*, \varrho, \varrho^*$
denote any sequence of scalars taken from $\fK$
  which satisfies (\ref{askwilb1}),
  (\ref{askwilb2}).
Applying
Lemma
\ref{paul12p2} we find
$P(\theta_{i-1},\theta_i)=0$ for $1 \leq i \leq d$,
where $P$ is from Definition
\ref{def:poly}.
By this and Lemma \ref{lem:pre}
 we find
 $\theta_0, \theta_1, \ldots, \theta_d$ is
$(\beta, \gamma)$-recurrent.
Therefore
 $\theta_0, \theta_1, \ldots, \theta_d$ is
$\beta$-recurrent.
By these remarks and since $d \geq 3$ we find
each of $\beta, \gamma, \varrho$ is uniquely
determined.
Interchanging $A,A^*$ in this argument
we find
 each of
$\gamma^*, \varrho^*$ is uniquely
determined. We have now proved the theorem for the case
$d\geq 3$.

\theorembreak Next assume $d\leq 2$.
 Let $\beta $ denote any scalar in $\fK$. If
$d=2$ define $\gamma=\theta_0-\beta \theta_1+\theta_2$ and if
$d\leq 1$ let $\gamma$ denote any scalar in $\fK$. If $d\geq 1$
define
\begin{eqnarray*}
\varrho =
 \theta_1^2-\beta\theta_1\theta_0
+\theta_0^2-\gamma (\theta_1+\theta_0)
\end{eqnarray*}
and if $d=0$ let $\varrho$ denote any scalar in $\fK$. Observe
$\theta_0, \theta_1, \ldots, \theta_d$ is $(\beta,
\gamma)$-recurrent; applying Lemma \ref{lem:pre}{\it (i)} we find
$P(\theta_{i-1},\theta_i)=0$ for $1 \leq i \leq d$. Applying Lemma
\ref{paul12p2} we find $\beta, \gamma, \varrho$ satisfy
(\ref{askwilb1}). Interchanging $A,A^*$ in the above argument we
find there exists scalars $\gamma^*, \varrho^*$ in $\fK$ such that
$\beta, \gamma^*, \varrho^*$ satisfy (\ref{askwilb2}). \qed

\noindent
We finish this section with a comment.

\begin{theorem}
{\rm \cite[Theorem 4.3]{qSerre}} \label{thm:lastcom}
 Let $d$
denote a nonnegative integer and let $V$ denote a vector space
over $\fK$ with dimension $d+1$. Let $(A,A^*)$ denote a Leonard
pair on $V$. Let the scalars $\theta_i,\theta^*_i$ be as in
Definition $\ref{primdef}$. Let
$\beta,\gamma,\gamma^*,\varrho,\varrho^*$ denote a sequence of
scalars taken from $\fK$ which satisfies $(\ref{askwilb1})$ and
$(\ref{askwilb2})$. Then the following (i)--(v) hold.
\begin{enumerate}
\item The expressions
\begin{equation}
\frac{\theta_{i-2}-\theta_{i+1}}{\theta_{i-1}-\theta_i},\qquad
\frac{\theta^*_{i-2}-\theta^*_{i+1}}{\theta^*_{i-1}-\theta^*_i}
\end{equation}
are both equal to $\beta+1$ for $2\le i\le d-1$.
\item $\gamma=\theta_{i-1}-\beta\theta_i+\theta_{i+1}, \qquad 1\le
i\le d-1$.
\item $\gamma^*=\theta^*_{i-1}-\beta\theta^*_i+\theta^*_{i+1},
\qquad
1\le i\le d-1$.
\item
$\varrho=
\theta^2_{i-1}-\beta\,\theta_{i-1}\,\theta_i
+\theta_{i}^2-\gamma\,(\theta_{i-1}+\theta_i),
\qquad 1 \leq i\leq d$.
\item
$\varrho^*=
\theta_{i-1}^{*2}-\beta\,\theta^*_{i-1}\,\theta^*_i
+\theta_{i}^{*2}-\gamma^*\,(\theta^*_{i-1}+\theta^*_i),
\qquad
1 \leq i\leq d$.
\end{enumerate}
\end{theorem}
\proof
{\it (iv), (v)} By
Corollary
\ref{lpscalars1} we have
$P(\theta_{i-1}, \theta_i)=0$
and
$P^*(\theta^*_{i-1}, \theta^*_i)=0$
for $1 \leq i \leq d$, where $P$ and $P^*$ are from
Definition
\ref{def:poly}.
\\
\noindent {\it (ii)} Combine Lemma \ref{lem:pre}{\it (ii)} with
part {\it (iv)} of this lemma.
\\
\noindent {\it (iii)} Similar to the proof of {\it (ii)} above.
\\
\noindent {\it (i)} The sequence $\theta_0, \theta_1, \ldots,\theta_d$
is $(\beta, \gamma)$-recurrent by {\it (ii)} so this sequence
is $\beta$-recurrent.
The sequence $\theta^*_0, \theta^*_1, \ldots,\theta^*_d$
is $(\beta, \gamma^*)$-recurrent by {\it (iii)} so this sequence
is $\beta$-recurrent. The result follows.
\qed

\section{The proof of Theorem 1.5}

In this section we prove Theorem
 \ref{lptheorem}.
We begin with an extension of
Lemma
\ref{paul12p2}.

\begin{lemma} \label{paulth14}
Let $\beta, \gamma, \varrho, \gamma^*, \omega, \eta$ denote
scalars in $\fK$. Then with reference to Definition
$\ref{primdef}$ and Definition $\ref{def:poly}$,
\begin{eqnarray} \label{thawnew}
A^2 A^*-\beta A A^*\!A+A^*\!A^2-\gamma\left( A A^*\!+\!A^*\!A
\right)-\varrho\,A^* &\equal& \gamma^*\!A^2+\omega A+\eta\,I
\end{eqnarray}
if and only if both
\begin{eqnarray}
\label{eq:cond1}
P(\theta_{i-1}, \theta_i)&=&0, \qquad 1 \leq i \leq d,
\\
\label{eq:cond2}
a^*_iP(\theta_{i}, \theta_i)&=&
\gamma^*\theta^2_i +\omega \theta_i + \eta,
 \qquad 0 \leq i \leq d.
\end{eqnarray}
\end{lemma}
\proof Let $L$ (resp. $R$) denote the expression on the left
(resp. right)
in
(\ref{thawnew}).
Observe
\begin{eqnarray}
L = \sum_{i=0}^d \sum_{j=0}^d E_iLE_j,
\qquad \qquad
R = \sum_{i=0}^d \sum_{j=0}^d E_iRE_j.
\label{RLmain}
\end{eqnarray}
Observe further that for $0 \leq i,j\leq d$ both
\begin{eqnarray}
E_iLE_j &=& P(\theta_i,\theta_j)E_iA^*E_j,
\label{eq:l}
\\
E_iRE_j &=& \delta_{ij}(
\gamma^*\theta^2_i +\omega \theta_i + \eta)E_i.
\label{eq:r}
\end{eqnarray}
First assume
(\ref{thawnew}), so that $L=R$.
We show
(\ref{eq:cond1}),
(\ref{eq:cond2}).
Concerning (\ref{eq:cond1}),
for $1 \leq i \leq d$ we have
$E_{i-1}RE_i=0$
by
(\ref{eq:r})
so
$E_{i-1}LE_i=0$.
By this and
(\ref{eq:l}) we find
$P(\theta_{i-1},\theta_i)E_{i-1}A^*E_i=0$.
Recall $E_{i-1}A^*E_i \not=0$
by Lemma
 \ref{leosystem}{\em (ii)}
so $P(\theta_{i-1},\theta_i)=0$. We now have (\ref{eq:cond1}).
Concerning (\ref{eq:cond2}), for $0 \leq i \leq d$ we have
$E_iLE_i=E_iRE_i$. Evaluating this using (\ref{eq:l}),
(\ref{eq:r}) and Lemma \ref{lem:ais}{\it (ii)}
 we find
$a^*_iP(\theta_i,\theta_i)=
\gamma^*\theta^2_i +\omega \theta_i + \eta $.
We now have
(\ref{eq:cond2}).
We have now proved the lemma in one direction.
To obtain the converse,
 assume
(\ref{eq:cond1}),
(\ref{eq:cond2}).
We show $L=R$.
By
(\ref{RLmain}), it suffices to
show
$E_iLE_j=E_iRE_j$ for
 $0 \leq i,j\leq d$.
For $0 \leq i \leq d$ we have
$E_iLE_i=a^*_iP(\theta_i,\theta_i)E_i$ by Lemma \ref{lem:ais}{\it
(ii)} and (\ref{eq:l}). Moreover
$E_iRE_i=(\gamma^*\theta^2_i+\omega \theta_i+\eta)E_i$ by
(\ref{eq:r}) so $E_iLE_i=E_iRE_i$ in view of (\ref{eq:cond2}). For
$1 \leq i \leq d$ we have $E_{i-1}LE_i=0$ by (\ref{eq:cond1}),
(\ref{eq:l}) and $E_{i-1}RE_i=0$ by
 (\ref{eq:r}) so
$E_{i-1}LE_i=
E_{i-1}RE_i$.
For $0 \leq i,j\leq d$ with $|i-j|>1$,
recall $E_iA^*E_j=0$ by Lemma
 \ref{leosystem}{\it (ii)} so
$E_iLE_j=0$ in view of
(\ref{eq:l}).
Also $E_iRE_j=0$ by
(\ref{eq:r})
so
 $E_iLE_j=
E_iRE_j$.
Apparently
 $E_iLE_j=
E_iRE_j$ for $0 \leq i,j\leq d$.
 By this and
(\ref{RLmain})
we find
$L=R$.
We now have
(\ref{thawnew}).
\qed

\begin{corollary} \label{lp2scalarsint}
Let
$\beta,\gamma,\gamma^*,\varrho, $
$\varrho^*,
\omega, \eta, \eta^*$
 denote
scalars in $\fK$. Then with reference to
 Definition $\ref{primdef}$
 and Definition $\ref{def:poly}$,
 the following (i), (ii) are equivalent.
\begin{enumerate}
\item The sequence
$\beta,\gamma,\gamma^*,\varrho,\varrho^*,
\omega, \eta, \eta^*$
satisfies {\rm (\ref{askwil1})} and {\rm (\ref{askwil2})}.
\item For $1 \leq i \leq d$ both
\begin{eqnarray*}
P(\theta_{i-1}, \theta_i)=0,
\qquad \qquad
P^*(\theta^*_{i-1}, \theta^*_i)=0
\end{eqnarray*}
and for $0 \leq i \leq d$ both
\begin{eqnarray}
\label{Pcond1}
a^*_iP(\theta_i, \theta_i)&=&
\gamma^*\theta^2_i +\omega \theta_i +\eta,
\\
\label{Pcond2}
 a_iP^*(\theta^*_i, \theta^*_i)&=&
 \gamma\theta^{*2}_i +\omega \theta^*_i +\eta^*.
\end{eqnarray}
\end{enumerate}
\end{corollary}
\proof
Apply Lemma
\ref{paulth14} to both $(A,A^*)$ and $(A^*,A)$.
\qed

\noindent {\bf Proof of Theorem
\ref{lptheorem}.}
Set $d+1=\mbox{dim} V$. For $d=0$ the result is trivial so assume
$d\geq 1$. For $\inint{i}{d}$, let $\theta_i$ and $E_i$ (resp.
$\theta^*_i$ and $E^*_i$) denote the eigenvalues and the primitive
idempotents of $A$ (resp. $A^*$) as in Definition \ref{primdef}.
Let $\mathcal D$ denote the $\fK$-subalgebra of $\mbox{End}(V)$
generated by $A$. By
 Theorem \ref{tdptheorem}
there exists scalars
$\beta, \gamma, \gamma^*, \varrho, \varrho^*$ in $\fK$
which satisfy
(\ref{askwilb1}), (\ref{askwilb2}).
We show there exist scalars
$\omega, \eta, \eta^*$ in $\fK$
such that
the sequence $\beta, \gamma, \gamma^*, \varrho, \varrho^*,
\omega, \eta, \eta^*$
 satisfies
 (\ref{askwil1}), (\ref{askwil2}).
Define
\begin{equation} \label{defrr}
U=A^2 A^*-\beta A A^*\!A+A^*\! A^2-\gamma\left( A
A^*\!+\!A^*\!A\right)-\varrho\,A^*.
\end{equation}
By
(\ref{askwilb1})
the element $U$ commutes with
$A$.
By this and since $A$ is multiplicity-free
we find
$U \in {\mathcal D}$.
Recall $\mathcal D$ has a basis $I,A, A^2, \ldots, A^d$
so
there exists a
polynomial
$f\in\fK[x]$
which has degree
at most $d$ and satisfies
$U=f(A)$.
Let $h$ denote the degree of $f$.
We show $h\leq 2$.
To do this we assume $h>2$ and get a contradiction.
We multiply each term in $U=f(A)$ on the left by $E^*_h$
and on the right by $E^*_0$.
We evaluate the result using
 (\ref{acomei}) and Lemma
\ref{eiarej}  to find $0=\alpha E^*_hA^hE^*_0$, where $\alpha$
denotes the leading coefficient of $f$. The scalar $\alpha$ is
nonzero by the construction and $E^*_hA^hE^*_0\not=0$ by Lemma
\ref{eiarej}{\it (ii)}. Therefore $0\not=\alpha E^*_hA^hE^*_0$ for
a contradiction. We have now shown $h\leq 2$, so there exist
scalars $\lambda, \omega, \eta $ in $\fK$ such that $U=\lambda
A^2+\omega A+\eta I$.

\theorembreak
 For the moment assume $d=1$. In this case $A^2$
is linearly dependent on $I, A$,
  so $\lambda, \omega, \eta $ can be chosen 
  such that $\lambda =\gamma^*$.
Next assume $d\geq 2$. We show
  $\lambda =\gamma^*$.
To do this, we
multiply each term in
$U=\lambda A^2+\omega A+\eta I$
on the left by $E^*_2$ and on the right by $E^*_0$.
We evaluate the result using
 (\ref{acomei}) and Lemma
\ref{eiarej}  to find
\begin{equation}
\label{eq:2a0}
\left( \theta^*_0-\beta\theta^*_1+\theta^*_2
\right)\,E^*_2\,A^2\,E^*_0 = \lambda \,E^*_2\,A^2\,E^*_0.
\end{equation}
Observe $E^*_2A^2E^*_0\not=0$ by Lemma \ref{eiarej}{\it (ii)}; by
this and (\ref{eq:2a0}) we find
$\theta^*_0-\beta\theta^*_1+\theta^*_2=\lambda$. Setting $i=1$ in
Theorem \ref{thm:lastcom}{\it (iii)} we find
$\theta^*_0-\beta\theta^*_1+\theta^*_2=\gamma^*$. We now see
 $\lambda=\gamma^*$.

\theorembreak
 Now $\beta, \gamma, \varrho, \gamma^*, \omega, \eta $ satisfy
(\ref{askwil1}). Interchanging the roles of $A$ and $A^*$ in the
argument so far, we find there exist scalars $\omega^*,\eta^*$ in
$\fK$ such that
\begin{equation} \label{askwil2a}
A^*{}^2\!A-\beta A^*\!AA^*\!+AA^*{}^2-\gamma^*\!\left(A^*\!A\!+\!A
A^*\right)-\varrho^*\!A = \gamma A^*{}^2+\omega^* A^*+\eta^*I.
\end{equation}
We show
$\omega^*=\omega$. We proceed as follows.
We first find the commutator of each
side of (\ref{askwil1}) with $A^*$. The result is
\begin{eqnarray*}
A^2A^*{}^2-\beta AA^*AA^*+\beta A^*AA^*A-A^*{}^2A^2
-\gamma(AA^*{}^2-A^*{}^2A)\\
=\gamma^*(A^2A^*-A^*A^2)+\omega (AA^*-A^*A).
\end{eqnarray*}
We next find the commutator
of each side of (\ref{askwil2a}) with $A$. The result is
\begin{eqnarray*}
A^*{}^2A^2-\beta A^*AA^*A+\beta AA^*AA^*-A^2A^*{}^2
-\gamma^*(A^*A^2-A^2A^*)\\
=\gamma(A^*{}^2A-AA^*{}^2)+\omega^* (A^*A-AA^*).
\end{eqnarray*}
Adding the last two equations and simplifying the result
we obtain
\begin{equation} \label{omegas}
0=(\omega-\omega^*)(AA^*-A^*A).
\end{equation}
Observe that $AA^*\not=A^*A$ since
\begin{equation} \label{fomegas}
E^*_1(AA^*-A^*A)E^*_0=(\theta^*_0-\theta^*_1)E^*_1 A E^*_0
\end{equation}
is nonzero. By (\ref{omegas}) and since $AA^*\not=A^*A$ we find
$\omega^*=\omega$.
Hence $\beta, \gamma^*, \varrho^*, \gamma, \omega, \eta^*$ satisfy
(\ref{askwil2}).

\theorembreak
 We have now shown
 there exists a sequence of scalars
$\beta, \gamma, \gamma^*, \varrho, \varrho^*, \omega, \eta,
\eta^*$ taken from $\fK$ which satisfies
 (\ref{askwil1}), (\ref{askwil2}).
We now assume $d\geq 3$ and show this sequence is
unique.
Let
$\beta, \gamma, \gamma^*, \varrho, \varrho^*,
\omega, \eta, \eta^*$ denote any
sequence of scalars in $\fK$ which satisfies
 (\ref{askwil1}), (\ref{askwil2}).
Then the sequence
$\beta, \gamma, \gamma^*, \varrho, \varrho^*$
satisfies
(\ref{askwilb1}), (\ref{askwilb2}) and is therefore
uniquely determined by Theorem
 \ref{tdptheorem}.
By (\ref{askwil1})
and since $I, A$ are linearly independent
we find
 $\omega, \eta$ are uniquely determined.
By (\ref{askwil2}) we find
 $\eta^*$ is uniquely determined.
We have now shown
the sequence
$\beta, \gamma, \gamma^*, \varrho, \varrho^*,
\omega, \eta, \eta^*$ is uniquely determined.
\qed

\noindent We finish this section with a comment. Let $V$ denote a
vector space over $\fK$ with finite positive dimension and let
$(A,A^*)$ denote a Leonard pair on $V$. Let
$\beta,\gamma,\gamma^*,\varrho,\varrho^*, \omega, \eta, \eta^*$
denote a sequence of scalars taken from $\fK$ which satisfies
(\ref{askwil1}), (\ref{askwil2}). Observe the scalars
$\beta,\gamma,\gamma^*,\varrho,\varrho^*$ satisfy {\it (i)--(v)}
in Theorem \ref{thm:lastcom}. Concerning $\omega, \eta, \eta^*$ we
have the following.

\begin{theorem} \label{extraeta}
Let $d$ denote a positive integer
and let $V$ denote a vector space over $\fK$ with dimension $d+1$.
Let $(A,A^*)$ denote a Leonard pair on $V$.
 Let
$\beta,\gamma,\gamma^*,\varrho, $  $\varrho^*,\omega, \eta,
\eta^*$ denote a sequence of scalars taken from $\fK$ which
satisfies $(\ref{askwil1})$, $(\ref{askwil2})$. Let the scalars
$\theta_i, \theta^*_i, a_i, a^*_i$  be as in
 Definition $\ref{primdef}$.
 For notational convenience, let $\theta_{-1}$ and
$\theta_{d+1}$ (resp.
 $\theta^*_{-1}$ and $\theta^*_{d+1}$)
denote scalars in $\fK$ which satisfy Theorem
$\ref{thm:lastcom}$(ii) (resp. Theorem $\ref{thm:lastcom}$(iii))
for $i=0$ and $i=d$. Then the following (i)--(iv) hold.
\begin{enumerate}
\item
$\omega =
a^*_i\,(\theta_i-\theta_{i+1})+a^*_{i-1}\,
(\theta_{i-1}-\theta_{i-2})-\gamma^*\,(\theta_i+\theta_{i-1}),
\qquad 1 \leq i \leq d$.
\item
$\omega = a_i\,(\theta^*_i-\theta^*_{i+1})+a_{i-1}\,
(\theta^*_{i-1}-\theta^*_{i-2})-\gamma\,(\theta^*_i+\theta^*_{i-1}),
\qquad 1 \leq i \leq d$.
\item
$\eta =
a^*_i\,(\theta_i-\theta_{i-1})\,(\theta_i-\theta_{i+1})
-\gamma^*\,\theta^2_i-\omega\,\theta_i,
\qquad  0 \leq i \leq d$.
\item
$\eta^*=
a_i\,(\theta^*_i-\theta^*_{i-1})\,(\theta^*_i-\theta^*_{i+1})
-\gamma\,\theta^{*2}_i-\omega\,\theta^*_i,
\qquad  0 \leq i \leq d$.
\end{enumerate}
\end{theorem}
\proof
{\it (iii)}
Let $i$ be given.
We claim
\begin{equation} \label{extraeq2}
P(\theta_i, \theta_i)=
(\theta_i-\theta_{i-1})(\theta_i-\theta_{i+1}),
\end{equation}
where $P$ is from
Definition
\ref{def:poly}.
To verify
(\ref{extraeq2}) for $1 \leq i \leq d$,
first eliminate $\theta_{i+1}$
using
$\theta_{i-1}-\beta \theta_i+\theta_{i+1}=\gamma$.
Evaluate the result using $P(\theta_i, \theta_i)=
(2-\beta)\theta^2_i-2\gamma \theta_i -\varrho$ and
$P(\theta_{i-1},\theta_i)=0$.
To verify
(\ref{extraeq2}) for $0 \leq i \leq d-1$,
first eliminate $\theta_{i-1}$
using
$\theta_{i-1}-\beta \theta_i+\theta_{i+1}=\gamma$.
Evaluate the result using $P(\theta_i, \theta_i)=
(2-\beta)\theta^2_i-2\gamma \theta_i -\varrho$ and
$P(\theta_{i},\theta_{i+1})=0$.
We now have
(\ref{extraeq2}).
Combining
(\ref{extraeq2})
with
(\ref{Pcond1})
we obtain the desired formula.
\\
\noindent {\it (iv)} Similar to the proof of {\it (iii)} above.
\\
\noindent {\it (i)}
Subtract {\it (iii)} (at $i$) from {\it (iii)} (at $i-1$) and simplify.
\\
\noindent {\it (ii)} Similar to the proof of {\it (i)} above.
\qed

\section{Concluding remarks}

We illustrate Theorem
\ref{lptheorem} by
computing the Askey-Wilson relations for the
Leonard pair in Example
\ref{krawtchouk}.

\begin{example} \label{krawtchouk2} \rm
Let $(A,A^*)$ denote the Leonard pair from Example
\ref{krawtchouk}. Referring to that example,
in the basis $y_0, y_1, \ldots, y_d$
the matrices for $A$ and $A^*$ have diagonal entries
\[
a_i=p(d-i)+(1-p)i, \qquad \theta_i^*=-i, \qquad  \qquad
 0\le i\le d.
\]
In the basis
 $y^*_0, y^*_1, \ldots, y^*_d$
the matrices for $A$ and $A^*$ have diagonal
entries
\[
\theta_i=i, \qquad a^*_i=-p(d-i)-(1-p)i,
\qquad
\qquad 0\le i\le d.
\]
\noindent
Define $\beta=2$,
$\gamma=\gamma^*=0$, $\varrho=\varrho^*=1$,
 $\omega=1-2p$, $\eta=pd$,
$\eta^*=-pd$. One readily verifies these scalars satisfy Corollary
\ref{lp2scalarsint}{\it (ii)}. Applying that corollary we find
\begin{eqnarray} 
A^2A^*-2AA^*A+A^*A^2 - A^*\! & \equal &
(1-2p)A+pdI,\\
A^{*2}A-2A^*AA^*+AA^{*2} - A & \equal &
(1-2p)A^*-pdI.
\end{eqnarray}
These are the Askey-Wilson relations for $(A,A^*)$.
\qed
\end{example}

\noindent
We conclude this paper with a kind of converse to Theorem
\ref{lptheorem}.
\begin{theorem}
Let $V$ denote a vector space over $\fK$ with finite positive
dimension. Let $A:V\to V$ and $A^*:V\to V$ denote linear
transformations. Suppose that:
\begin{itemize}
\item There exists a sequence of scalars
$\beta,\gamma,\gamma^*,\varrho,\varrho^*,\omega,\eta,\eta^*$ taken
from $\fK$ which satisfies {\rm (\ref{askwil1})}, {\rm
(\ref{askwil2})}.
\item $q$ is not a root of unity, where $q+q^{-1}=\beta$.
\item Each of $A$ and $A^*$ is multiplicity-free.
\item There does not exist a subspace $W\subseteq V$ such that
$W\not=0$,
$W\not=V$,
$AW\subseteq W$,
$A^*W\subseteq W$.
\end{itemize}
Then $(A,A^*)$ is a Leonard pair on $V$.
\end{theorem}
\proof By
\cite[Theorem 3.10]{qSerre} the pair $(A,A^*)$ is
a tridiagonal pair on $V$ in the sense of \cite{TD00}.
Now by
\cite[Lemma 2.2]{qSerre} and since
each of $A, A^*$
is multiplicity-free, we find
$(A,A^*)$ is a Leonard pair on $V$.
 \qed



\noindent Paul Terwilliger \hfil\break
\noindent Department of Mathematics \hfil\break
\noindent University of Wisconsin \hfil\break
\noindent 480 Lincoln Drive \hfil\break
\noindent Madison, WI 53706 USA \hfil\break
\noindent email: terwilli@math.wisc.edu \hfil\break

\medskip
\noindent Raimundas Vidunas \hfil\break
 \noindent Dept. Mathematics and Computer Science \hfil\break
 \noindent RUCA, Antwerp University \hfil\break \noindent Middelheimlaan 1
\hfil\break \noindent 2020 Antwerp, Belgium \hfil\break \noindent
email: Raimundas.Vidunas@ua.ac.be  \hfil\break

\end{document}